\documentclass[titlepage,a4paper]{article}
\usepackage[T1]{fontenc}
\usepackage[latin1]{inputenc}
\usepackage{amsmath}
\usepackage{amssymb}
\usepackage{amsthm}

\newtheorem{Def}{Definition}[section]

\newtheorem{rem}{Remark}[section]
\newtheorem{Th}{Theorem}[section]
\newtheorem{Aux}{Auxiliary Lemma}[section]
\newtheorem{lem}{Lemma}[section]

\newcommand{\RR}{\mathbb{R}}
\newcommand{\EE}{\mathbb{E}}
\newcommand{\PP}{\mathbb{P}}
\newcommand{\MM}{\mathbb{M}}
\newcommand{\QQ}{\mathbb{Q}}
\newcommand{\NN}{\mathbb{N}}
\newcommand{\TT}{\mathbb{T}}

\newcommand{\Cov}{\mathrm{Cov}}
\newcommand{\cB}{{\cal B}}
\newcommand{\cA}{{\cal A}}
\newcommand{\cP}{{\cal P}}
\newcommand{\cL}{{\cal L}}
\newcommand{\cC}{{\cal C}}
\newcommand{\cE}{{\cal E}}
\newcommand{\cH}{{\cal H}}

\newcommand{\cF}{{\cal F}}

\newcommand{\cG}{{\cal G}}

\newcommand{{{\cadlag}}}{c\`adl\`ag}

\newcommand{\aN}{{{^*}\mathbb{N}}}

\renewcommand{\baselinestretch}{2}

\begin{document}

\pagestyle{plain}

\title{The fairest price of an asset in an environment of temporary arbitrage }

\author{Frederik S. Herzberg\thanks{Mathematical Institute, University of Oxford, Oxford OX1 3LB, United Kingdom, Telephone (reception): +44 (0)1865 273525, Facsimile: +44 (0)1865 273583, e-mail: {\tt herzberg@maths.ox.ac.uk}. The author gratefully acknowledges funding from the German National Academic Foundation ({\em Studienstiftung des deutschen Volkes}) and a PhD grant of the German Academic Exchange Service ({\em Doktorandenstipendium des Deutschen Akademischen Austauschdienstes}).}}

\renewcommand{\baselinestretch}{1} \maketitle

\maketitle

\begin{abstract} In practice there are temporary arbitrage opportunities arising from the fact that prices for a given asset at different stock exchanges are not instantaneously the same. We will show that even in such an environment there exists a ``fairest measure'' (instead of a martingale measure), albeit not necessarily unique. For this end, we define and analyse quantitative notions of unfairness in complete as well as incomplete market settings.\end{abstract}

\renewcommand{\baselinestretch}{2}

\section{Introduction}

Assets that are traded simultaneously at several different stock exchanges, such as currencies, attain different prices at the very same moment. This entails a short-term arbitrage opportunity which can be exploited with positive probability. On the other hand, there is of course a very strong correlation among the asset price vectors of different stock exchanges. 

In this framework of local arbitrage, we are going to show the existence of a measure that, although not being a martingale measure for the asset prices, minimises the ``distance'' of the asset price vectors from being a martingale.

Although we will not construct this ``fairest measure'' on the original probability space $\left(\Gamma,(\cG_t)_{t\in[0,1]},\PP\right)$, we will show that there is such a measure on an adapted probability space $\left(\Omega,(\cF_t)_{t\in[0,1]},\QQ\right)$ which is {\em equivalent} \cite{FK02,K88} to $\Gamma$. Informally this equivalence relation means that in some sense every reasonable probabilistic assertion that holds on $\left(\Gamma,(\cG_t)_{t\in[0,1]},\PP\right)$ is also true on $\left(\Omega,(\cF_t)_{t\in[0,1]},\QQ\right)$ (and {\em vice versa}).

To the knowledge of the author, no attempt has been made so far to develop quantitative notions of unfairness for price processes. The theory of fairness functionals for surfaces \cite{W,KR} is not suitable for our purposes in a straightforward manner, since our price processes will have rough paths and therefore do not even allow to apply the surface notion of fairness pathwise. However, it is well possible to generalise the approach taken for the definition of complete-market unfairness measures in this paper and apply, given a discounted price process $z$ with natural filtration $\cF$ and a starting time $t_0\geq 0$, fairness functionals pathwise to the process $t\mapsto \EE\left[\left.z_{t+t_0}\right|\cF_{t_0}\right]$. In order to make this definition independent from the choice of a starting point, $t_0$, it is necessary to integrate over $t_0$. That said, if we require our measures of unfairness to be $1$-homogneneous (as we do in the complete market case) or even to remain unchanged under multiplication with constants (as is the case for our notion of unfairness in incomplete markets), the theory of fair surfaces does not give us many new suggestions for natural methods of quantifying unfairness. 

In any case, this paper aims at a rigorous generalisation of the commonly accepted view that fair discounted price processes should be martingales (as they do not admit arbitrage) and {\em vice versa}, thus we will not be able to refrain from an irreducibly probabilistic definition of unfairness.

Our proofs will employ techniques from nonstandard analysis; these have already been successfully applied to mathematical finance in the work of Cutland, Kopp and Willinger \cite{CKW93}.

\section{Definitions and main results}

\begin{Th}\label{th1} Let $n\in\NN$, $p>0$ and $c,N>0$. Consider an adapted probability space $\left(\Gamma,\cG, \PP\right)$ and $\cB[0,1]\otimes\cG_1$-measurable geometric Brownian motions with constant multiplicative drift $\tilde g_i:\Gamma\times[0,1]\rightarrow \RR^{d\cdot n}$, $i\in\{1,\dots,n\}$. Then there exist processes $g_i$, $i\in\{1,\dots,n\}$ on an adapted probability space $\left(\Omega,\cF,\QQ\right)$, equivalent to the processes $\tilde g_i$ on $\Gamma$ (in the sense of adapted equivalence \cite{FK02}), such that there is a probability measure $\MM_m$ on $\Omega$ minimising $$Q\mapsto m_\Omega(Q, g):=\sum_{i=1}^n \int_0^1\int_s^1 \int_\Omega \left| \left(g_i\right)_s - \EE_Q\left[\left. \left(g_i\right)_t\right| \cF_s\right] \right|^pdQ\ dt \ ds$$ in the class of measures $$\cC(\Omega,g):=\left\{Q:\cF_1\rightarrow[0,1] \ : \ \begin{array}{c} Q\text{ probability measure}, \\ \forall A\in\cF_1 \quad\frac{1}{N}\cdot \QQ(A) \leq {Q(A)}\leq N \cdot\QQ(A),\\ \forall i\neq j\in\{1,\dots, n\}\\ \int_0^1 \frac{\Cov_Q\left(\left(g_i\right)_s,\left(g_j\right)_s\right) }{\EE_Q\left|\left(g_ig_j\right)_s\right|}ds\geq c \end{array} \right\}.$$ 
Analogously, there is a probability measure $\MM_n$ minimising $$Q\mapsto n_\Omega\left(Q,g\right)=\int_0^1\EE_\QQ\left|\left.\frac{1}{g_t}\frac{d}{du}\right|_{u=0}\EE\left[g_{t+u}\left.\right|\cG_t\right]\right|dt$$ (where $n_\Omega(Q,g)$ is defined to be $+\infty$ if the derivative in $0$ in this definition does not a.s. exist as a continuous function in $t$) in the class $\cC(\Omega,g)$. Moreover, $$\inf_{\cC(\Omega,g)}m_\Omega(\cdot, g)\leq \inf_{\cC\left(\Gamma,\tilde g\right)}m_\Gamma\left(\cdot,\tilde g\right)$$ as well as $$\inf_{\cC(\Omega,g)}n_\Omega(\cdot, g)\leq \inf_{\cC\left(\Gamma,\tilde g\right)}n_\Gamma\left(\cdot,\tilde g\right).$$

\end{Th}

The proof for this Theorem can be split into the following Lemmas which might also be interesting in their own right:

\begin{lem} \label{lem1}Using the notation of the previous Theorem \ref{th1}, for any hyperfinite adapted space $\Omega$ \cite{K88}, $$\inf_{\cC(\Omega,g)}m_\Omega(\cdot, g)\leq\inf_{\cC\left(\Gamma,\tilde g\right)}m_\Gamma\left(\cdot,\tilde g\right).$$
\end{lem}

\begin{lem} \label{lem2}Under the assumptions of Theorem \ref{th1} and choosing $\Omega$ to be any hyperfinite adapted space, the infimum of $m_\Omega (\cdot, g)$ on $\cC(\Omega, g)$ is attained by some measure $\MM_m\in \cC(\Omega,g)$.
\end{lem}

Although they look very similar, it is technically slightly more demanding to prove the following to Lemmas (which in turn obviously entail the second half of the Theorem, i.e. the assertions concerned with the map $n$):

\begin{lem} \label{lem3}Using the notation of the previous Theorem \ref{th1}, for any hyperfinite adapted space $\Omega$ \cite{K88}, $$\inf_{\cC(\Omega,g)}n_\Omega(\cdot, g)\leq\inf_{\cC\left(\Gamma,\tilde g\right)}n_\Gamma\left(\cdot,\tilde g\right).$$
\end{lem}

\begin{lem} \label{lem4}Under the assumptions of Theorem \ref{th1} and choosing $\Omega$ to be any hyperfinite adapted space, the infimum of $n_\Omega (\cdot, g)$ on $\cC(\Omega, g)$ is attained by some measure $\MM_n\in \cC(\Omega,g)$.
\end{lem}

Easy results are

\begin{lem} A semimartingale $x$ is a $P$-martingale on $\Omega$ if and only if $m_\Omega (P,x)=0$. The function $m_\Omega(P,\cdot)^{\frac{1}{p}}$ on the space of measurable processes of $\Omega$ satisfies the triangle inequality and is $1$-homogeneous. For $p=2$, it defines an inner product on the space $$ \cE:=\left\{x:\Omega\times[0,1]\rightarrow \RR^d \ : \ x\text{ measurable},\quad m(P,x)<+\infty\right\}$$ which becomes a Hilbert space by this construction.
\end{lem} 

\begin{lem} A semimartingale $x$ is a $P$-martingale on $\Omega$ if and only if $n_\Omega (P,x)=0$. The function $m_\Omega(P,\cdot)^{\frac{1}{p}}$ on the space of measurable processes of $\Omega$ remains unchanged when multiplying the argument by constants ({\em scale invariance}).
\end{lem}

We can generalise this to a 

\begin{Def} Let $\Omega$ be an adapted probability space. Define $$\cL\left(\Omega,\RR^d\right):= \left\{x:\Omega\times [0,1]\rightarrow\RR^d \ : \ x \text{ measurable}\right\}.$$ A function $\Upsilon:\cL\left(\Omega,\RR^d\right)\rightarrow[0,+\infty)$ is an {\em incomplete market notion of unfairness} if and only if it satisfies the triangle inequality, is $1$-homogeneous, and assigns $0$ to a semimartingale $y$ if and only if $y$ is a martingale.
$\Upsilon$ is said to be a {\em complete market notion of unfairness} if and only if it remains unchanged under multiplication by constants and $\Upsilon$ vanishes exactly for those semimartingales that are in fact maritngales.
\end{Def}

\begin{rem} The distinction between complete and incomplete markets can be justified by the following reasoning: If it is, under assumption of completeness, possible to buy as much of an asset as on eintends to, multiplication by a constant does not enhance the arbitrage opportunity at all; therefore, for complete markets, a suitable notion of unfairness should be scale invariant in that it does not change under multiplication of the argument by constants. Conversely, if we work in an incomplete market, there are restrictions on the aquisition of assets, therefore one has by genericity to assume that one will only be concerned with the development of one stock's price. Thus, a $1$-homogeneous notion of unfairness is most suitable in this situation.
\end{rem}

\section{Proofs}

\begin{proof}[Proof of Lemma \ref{lem1}] One applies Fajardo's and Keisler's model theory of stochastic processes \cite{K88,FK02}. It is sufficient to find an $\RR^{d\cdot n}$-valued process $g$ such that there is for any $P\in\cC\left(\Gamma,\tilde g\right)$ a $Q\in\cC(\Omega,g)$ satisfying $m_\Gamma\left(P,\tilde g\right)=m_\Omega(Q,g)$. Due to the universality of hyperfinite adapted probability spaces \cite[Adapted Universality Theorem]{FK02}, there is a process $g$ on $\Omega$ with the property that $\left(\Gamma,\tilde g\right)\equiv (\Omega,g)$ (where $\equiv$ denotes adapted equivalence). Now consider any $P\in\cC\left(\Gamma,\tilde g\right)$. Due to the (even uniformised) $\PP$-continuity of all elements of $\cC\left(\Gamma,\tilde g\right)$, there exists -- due to the Radon-Nikod\'ym Theorem -- a $\PP$-density $\frac{dP}{d\PP}$ for $P$. Now, $\Omega$ is saturated as an adapted probability space \cite[Saturation Theorem 3A.3]{FK02}, hence there exists a function $f:\Omega\rightarrow \RR$ such that $$\left(\Gamma,\left(\tilde g, \frac{dP}{d\PP}\right)\right)\equiv \left(\Omega,(g,f)\right).$$ 
Hence the covariances of the components of $g$ and $\tilde g$ coincide and the finite-dimensional distribution of $\tilde g$ under $P$ is the same as the finite-dimensional distribution of $g$ under $Q$. Thus $$\left(\tilde g, \Gamma, P\right)\equiv_0 \left( g,\Omega, Q \right)$$ in the notation of Hoover and Keisler \cite{HK84} or Fajardo and Keisler \cite{FK02}. Now, according to Girsanov's Theorem, $\ln g_i$ and $\ln \tilde g_i$ are Brownian motions with constant drift with respect to $P$ and $Q$, respectively, for any $i\in\{1,\dots,n\}$. Hence, all the $g_i,\tilde g_i$ are Markov processes with respect to the measures $P$ and $Q$, respectively, as well. However, due to Hoover and Keisler \cite{HK84}, $\equiv_0$ and $\equiv $ coincide for Markov processes and we do already know $(g,\Omega,Q)\equiv_0\left(\tilde g,\Gamma, Q\right)$. Thus, $$(g,\Omega,Q)\equiv\left(\tilde g,\Gamma, Q\right)$$ which finally yields $$m_\Omega(Q,g)=m_\Gamma\left( P,\tilde g\right).$$ This suffices to prove the Lemma.

\end{proof}

\begin{proof}[Proof of Lemma \ref{lem2}] Use the Transfer Principle and a lifting argument. For finite probability spaces and finite time lines, the desired probability measure is the solution of a finite-dimensional optimisation problem on a closed subset of the respective unit cube for a continuous function; hence the existence is clear. Transfer yields the generalisation to hyperfinite adapted spaces and the existence of optimal internal probability measures $Q$ minimising $M_\Omega(\cdot,G)$ where $G$ is internal and $M$ is the straightforward internal analogue of $M$. To exploit that result for our case, one needs to consider liftings $G$ for $g$ and $F$ for the $\QQ$-density $f$ of any measure $Q\in\cC(\Omega,g)$. Then one has to achieve \begin{equation} \label{lift} °m_\Omega\left(Fd\mu,G\right)=m_\Omega\left(fd\QQ,g\right)\end{equation} ($\mu$ denoting the internal normalised counting measure on $\Omega$, i.e. $L(\mu)=\QQ$) via the Adapted Lifting Theorem \cite{FK02} and the following 
\begin{Aux} \label{aux} For all probability spaces $(\Omega, \cG,\QQ)$, random variables $x$, $\QQ$-probability densities $f$ and $\sigma$-subalgebras $\cH$ of $\cG$, $$\EE_{fd\QQ}\left[x|\cH\right]=\frac{\EE_\QQ\left[fx|\cH\right]}{\EE_\QQ\left[f|\cH\right]}$$ $\QQ$-almost surely.
\end{Aux}
\begin{proof}[Proof of the Auxiliary Lemma \ref{aux}] First of all, we need to show that -- after choosing a representative of the conditional $\QQ$-expectation of $f$ with repect to $\cH$ -- $A:=\left\{\EE_\QQ\left[f|\cH\right]= 0\right\}$ is a $\QQ$-null set. Observe however, that $A$ is a $\QQ$-null set if and only if $A$ is $fd\QQ$-null set, since $fd\QQ\in \cC(\Omega, g)$ (and thus $fd\QQ$ and $\QQ$ are -- even uniformly -- equivalent). Now, \begin{eqnarray*}(fd\QQ)(A)&=&\int_Afd\QQ =\int_A\EE_\QQ\left[f|\cH\right]d\QQ \\ &=& \int_{\left\{\EE_\QQ\left[f|\cH\right]= 0\right\}}\EE_\QQ[f|\cH]d\QQ=0.\end{eqnarray*} The rest of the proof is as easy a calculation: For all $B\in\cH$, \begin{eqnarray*} \int_B\frac{\EE_\QQ[fx|\cH]}{\EE_\QQ[f|\cH]}\cdot fd\QQ &=& \int_B\EE_\QQ\left[\left.f \cdot\frac{\EE_\QQ[fx|\cH]}{\EE_\QQ[f|\cH]}\right|\cH\right]d\QQ \\ &=& \int_B\frac{\EE_\QQ[fx|\cH]}{\EE_\QQ[f|\cH]}\cdot \EE_\QQ[f|\cH]d\QQ =\int_B{\EE_\QQ[fx|\cH]}d\QQ \\ &=&\int_B fxd\QQ\end{eqnarray*}
\end{proof}
Now for the proof of equation (\ref{lift}), consider any $i\in\{1,\dots,n\}$ and let $\cP_s$ be the internal partition generating the internal algebra $\cA_s$ for any $s\in \TT$, $\TT$ being the hyperfinite time line such that $\Omega={\Omega_0}^\TT$ for some internal $\Omega_0$. Then, due to the shape of internal conditional expectations, we have \begin{eqnarray*}&& °m_\Omega(Fd\mu,G)\\&=& °\int_\TT \int_{\TT\cap\{\cdot\geq s\}}\int_\Omega \left| (G_i)_s - \frac{\sum_{A\in\cP_s}\chi_A\sum_{\omega'\in A}(G_i)_t(\omega')F(\omega')\mu\{\omega'\}}{\sum_{B\in\cP_s}\chi_B\sum_{\omega''\in B}F(\omega'')\mu\{\omega''\}}\right|^p\\ && Fd\mu \ dt \ ds\end{eqnarray*} and we observe that \begin{eqnarray*} &&\left(\sum_{C\in\cP_s}\chi_C\sum_{\omega'''\in C}\mu\{\omega'''\}\right)\sum_{D\in\cP_s}\chi_D\frac{\sum_{\omega^{(4)}\in D}F\left(\omega^{(4)}\right)\mu\left\{\omega^{(4)}\right\}}{\sum_{\omega^{(5)}\in D}\mu\left\{\omega^{(5)}\right\}} \\ &=&\sum_{B\in\cP_s}\chi_B\sum_{\omega''\in B}F(\omega'')\mu\{\omega''\}.\end{eqnarray*} 
Therefore, using the well-known Loeb integration theory as developed by Loeb\cite{L75} and Anderson \cite{A76}, as well as the previous Auxiliary Lemma \ref{aux}, \begin{eqnarray*}&& °m_\Omega(Fd\mu,G)\\&=& °\int_\TT \int_{\TT\cap\{\cdot\geq s\}} \int_\Omega \left| (G_i)_s - \frac{\EE_\mu\left[\left.F\cdot (G_i)_t\right|\cA_s\right] }{\EE_\mu\left[F|\cA_s\right]}\right|^p Fd\mu \ dt \ ds \\ &=& \int_0^1 \int_s^1 \int_\Omega \left| (g_i)_s - \frac{\EE_{L(\mu)}\left[\left. f\cdot (g_i)_t\right|\cF_s\right] }{\EE_{L(\mu)}\left[f|\cF_s\right]}\right|^p fdL(\mu) \ dt \ ds \\ &=& \int_0^1 \int_s^1 \int_\Omega \left| (g_i)_s - \frac{\EE_{\QQ}\left[\left. f\cdot (g_i)_t\right|\cF_s\right] }{\EE_{\QQ}\left[f|\cF_s\right]}\right|^p dQ \ dt \ ds \\ &=& \int_0^1 \int_s^1 \int_\Omega \left| (g_i)_s - {\EE_{fd\QQ}\left[\left. (g_i)_t\right|\cF_s\right] } \right|^p dQ \ dt \ ds \\ &=& \int_0^1 \int_s^1 \int_\Omega \left| (g_i)_s - {\EE_{dQ}\left[\left. (g_i)_t\right|\cF_s\right] } \right|^p dQ \ dt \ ds .\end{eqnarray*} This proves equation (\ref{lift}).
\end{proof}

For the proof of Lemmas \ref{lem3} and \ref{lem4} we will have to use the following Auxiliary Lemma that is necessary because there is nothing like an immediate internal analogue $N$ for the map $n$ as it was the case with $M$ and $m$.

\begin{Aux}\label{aux2} Suppose $\Omega={\Omega_0}^{\TT_H}$ where $\TT$ is the hyperfinite time line of cardinality $H!+1$. Then there is an $h\in\aN\setminus\NN$ such that $h\leq H$ and \begin{eqnarray*} &&L(\mu)\text{-a.s.}\quad \forall s\in\left\{0,\frac{1}{h!},\dots,1-\frac{1}{h!}\right\}\forall i\in\{1,\dots,n\}\subset \NN\\ && °\left(\frac{\EE\left[\left. F\cdot \left(G_i\right)_{s+\frac{1}{h!}}\right|\cA_s\right] \frac{1}{(G_i)_s}-\EE\left[F|\cA_s\right]}{\frac{1}{h!}\EE[F|\cA_s]}\right)\\ &=&\left.\frac{d}{du}\right|_{u=0} \frac{1}{(g_i)_{°s}\EE\left[f|\cF_{°s}\right]}\EE\left[\left. f\cdot (g_i)_{°s+u}\right|\cF_{°s}\right]\end{eqnarray*}
\end{Aux}
\begin{proof}[Proof of the Auxiliary Lemma \ref{aux2}] The proof is a saturation argument which we can exploit by using the fact that by definition of $n$, $u\mapsto \EE_{fd\QQ}[(g_i)_{t+u}|\cF_t]$ must be a.s. differentiable in $0$ for all $t\in[0,1]$ and $i\in\{1,\dots,n\}$ with the derivative behaving continuously in $t$ if $fd\QQ$ should be any candidate for minimising $n$. By virtue of the Auxiliary Lemma \ref{aux} we can deduce from this fact that $u\mapsto \EE_\QQ\left[f\cdot(g_i)_{t+u}|\cF_t\right]$ must be a.s. differentiable as well -- and the derivative will again behave continuously in $t$, since $\left(\EE[z|\cG_t]\right)_{t\in[0,1]}$ is a continuous martingale for any $\cG_1$-measurable random variable $z$, and any filtration $\cG$ coming from a continous martingale itself.
\end{proof}

With this Auxiliary Lemma \ref{aux2}, the proofs of Lemmas \ref{lem3} and \ref{lem4} can be carried out exactly in the same manner as we have proven Lemmas \ref{lem1} and \ref{lem2}.

Also, the uniform equivalence condition implies that we may assume that the fairest internal measure $\nu$ on the hyperfinite adapted space does not assign $0$ to any of the elementary events $\{\omega\}$ for $\omega\in\Omega={\Omega_0}^\TT$. Therefore, on this measure the local minimum of $m$ is attained for finite (bounded) $|\TT|$. Thus we obtain first order conditions on the density of $Q=L(\nu)$ with respect to $\QQ=L(\mu)$ for the case of finite $\TT$ which we can then generalise to the hyperfinite case.
\section{Final remarks}

Albeverio and Steblovskaya \cite{AS02, AS02a} have introduced models of financial markets in which the asset prices are not independent from each other but satisfy -- up to a Brownian motion -- some linear equation. In this framework, after imposing certain assumptions on the coefficients in the governing equations of the market, the existence of a martingale measure and hence the existence of (pseudo-)prices could be shown. However, in general there is no reason to expect the (empirical) prices of an asset traded at several stock exchanges to satisfy such a linear relation -- even up to a Brownian motion. Especially, we intend to deal with cases where there is, indeed, arbitrage. We have got no advantage as far as the uniqueness of the fairest measure is concerened. Neither in the setting of Albeverio and Steblovskaya nor in ours can uniqueness be expected in general.

\end{document}